\documentclass[12pt]{article}
\usepackage{geometry}                
\usepackage[all]{xy}
\usepackage{graphicx}
\usepackage{amssymb}
\usepackage{epstopdf}
\begin{document}
 \newcommand{\bfz}{\textbf{0}}
  \newcommand{\bfu}{\textbf{1}}
  \newcommand{\mR}{\mathcal{R}}
   \newcommand{\mI}{\mathcal{I}}
    \newcommand{\I}{\mathcal{I}}
 \newtheorem{definition}{Definition}
\newtheorem{lemma}{Lemma}
\newtheorem{theorem}{Theorem}
\newtheorem{remark}{Remark}
\newtheorem{example}{Example}
\author{Vittorio Cafagna and Gianluca Caterina}
  \title{Spectra of fuzzy sets of $S^1$}
    \maketitle

\abstract{In this article we explore some aspects of Fourier analysis on fuzzy sets. This work was part of a broader project jointly with Professor Vittorio Cafagna which spanned over a few years between 1998 and 2005. Professor Cafagna passed away unexpectedly in 2007. His intellectual breadth and inspiring passion for mathematics is still very well alive today.}
  \section{Introduction}
 Building on several papers from  the mid-thirties to beginning of fifties, from specialists of mathematical analysis like Verblunsky \cite{Ver}, Szasz \cite{Szasz1}, \cite{Szasz2}, Friedman \cite{Fried}, Geronimus \cite{Gero}, Ghizzetti \cite{Ghizz}, we give a complete characterization of the Fourier coefficients of fuzzy subsets of the unit circle. The main idea is to define a sort of {\sl nonlinear Fourier transform} which associates to a fuzzy subset of the unit circle the sequence of the Taylor coefficients of a suitable holomorphic function on the unit disk. It turns out that these {\sl nonlinear coefficients} have a very nice description as  sequences of complex numbers verifying a sequence of semi-algebraic conditions. The original Fourier coefficients are related to the nonlinear coefficients by a sequence of algebraic maps. The key role in this characterization is played by the nonlinear coefficients of a special family of crisp subsets: the characteristic functions of union of  $n$ arcs in the unit circle, called the {\sl crisp subsets of order n}. Highlighting a procedure which is central for the proof, we obtain that one can build, for every fuzzy subset $f$ of the unit circle, a sequence $\{\chi_n\}_{n\in\mathbb N}$ of crisp subsets of order $n$ such that the first $n$ Fourier coefficients of $\chi_n$ are the same as the first $n$ Fourier coefficients of $f$. The implications for fuzzy theory seem quite surprising: {\sl the spectrum of a fuzzy set is totally determined by the spectra of a sequence of crisp sets} (of course within the assumptions of the theorem: for the moment, summable fuzzy subsets of the unit circle.) We hope that this result may contribute to the vivid debate, initiated by the well-known paper by Elkan (\cite{Elk1}), {\sl The Paradoxical Success of Fuzzy Logic} (see also  \cite{Elk2} for different critics from the scientific community to the author's thesis). The main point of Elkan is that, despite the highly successful  applications of fuzzy logic to heuristic control, the conceptual machinery employed collapses to boolean logic.  We do not intend, and certainly we do not have the skills, to take a position in the {\sl pol\'emique}. But we hope that our result might be of some interest for the defenders of the one or the other position. Another issue, also highly debated in fuzzy theory is that of {\sl defuzzification}. We express the hope that the  {\sl spectral approximation} given by our sequences of crisp subsets of order $n$ might possibly suggest some new points of view on the concept of defuzzification. Finally, there is a research area which seems to be more than superficially related  to the argument of this chapter: spectral logics (see, e.g., the recent book \cite{moraga}). Even if, in that case, the main issue is to suggest a description of many-valued logics in terms of orthogonal wavelets bases (very far, both in content and methods, from our nonlinear Fourier technique), it is undeniable that we do share with the proponents of this research project the common language of harmonic analysis and the hope that a {\sl spectral} representation of logic might be of some interest.

\section{Preliminary notions of Fourier Analysis}
We review briefly the basic notions of Fourier analysis which we will need in the course of the chapter. We refer to \cite{D-McK} and \cite{Foll1} for the whole framework, missing parts of proofs only sketched and  concepts only alluded to.   All we need of functional analysis (and, of course, a lot more) can be found in the classical treatise of Dunford and Schwartz \cite{DS}. Some {\it Kindergarten} ideas of Fourier analysis on locally compact abelian groups are occasionally used here and there. Good references are \cite{Loomis}, \cite{Foll2} and again \cite{D-McK}.

Let $S^1=\{z\in\mathbb C:\parallel z\parallel=1\}$ be the unit circle, parametrized by $t\in[0,2\pi)$. Define  $L^1(S^1,\mathbb C)=\{f:S^1\to\mathbb C: \int_0^{2\pi} |f(t)| dt<\infty\}$, which is a Banach space with respect to the norm $\parallel f\parallel_1={1\over 2\pi}\int_0^{2\pi} |f(t)| dt$. Functions in $L^1(S^1,\mathbb C)$ will be referred to as {\sl  summable}. To every summable function on $S^1$ we can associate the sequences of its Fourier coefficients

$$c_k(f)={1\over 2\pi}\int_0^{2\pi} f(t)e^{ikt}dt,\,\,\,k\in\mathbb Z$$

$$a_0(f)={1\over \pi}\int_0^{2\pi} f(t) dt,\,\,\,a_k={1\over \pi}\int_0^{2\pi} f(t)coskt dt,\,\,\,k\in\mathbb N^*$$

$$b_k(f)={1\over \pi}\int_0^{2\pi} f(t)sinkt dt,\,\,\,k\in\mathbb N^*$$

\noindent  We will usually take the freedom of writing simply $c_k$ instead of $c_k(f)$, when then are no doubts of which function we are referring to. Euler's formula ($e^z=cosz+isinz$) gives the relation $c_k={a_k+ib_k\over 2}$. Remark that in the particular case that $f$ is real-valued, the sequence $\{c_k\}_{k\in\mathbb Z}$ is hermitian, i.e. $c_{-k}=\bar c_k$. The (for now) formal series $\sum_{k\in\mathbb Z}c_ke^{ikt},\,\, \sum_{k\in\mathbb N}a_kcoskt$ and $\sum_{k\in\mathbb N^*}b_ksinkt$ are called, respectively, the {\sl Fouries series}, the {\sl cosine Fourier series} and the {\sl sine Fourier series} of $f$.

The fundamental theorem of the theory is
\begin{theorem}{{\bf (Fourier's Theorem)}} Every summable function on $S^1$ is the sum of its Fourier series:
\end{theorem}

$$f(t)=\sum_{k\in\mathbb Z}c_ke^{ikt}$$
\begin{flushright}$\square$\end{flushright}

\begin{definition} 
If we denote by $\mathbb C^\mathbb Z$ the space of all sequences  $\{\alpha:\mathbb Z\to\mathbb C\}$, the map $f\mapsto \hat{f}$ from $L^1(S^1,\mathbb C)$ to $\mathbb C^\mathbb Z$ defined by

$$\hat{f}(k)=c_k(f)$$  

\noindent is called the {\sl Fourier transform}
and the sequence $\hat f$  {\sl the spectrum} of $f$.
\end{definition}

\begin{remark}{\bf (A glimpse of AHA: Abstract Harmonic Analysis)} Take a locally compact abelian group $G$ and define the set $\hat{G}=\{f:G\to S^1\}$, f a continuous group homomorphism (here the group operation on $S^1$ is just the multiplication of complex numbers). An element of $\hat{G}$ is called a {\sl character} of $G$. It can be shown (see, e.g. \cite{Loomis} or \cite{Foll2}) that the set of all characters of $G$ is itself a locally compact abelian group, called the {\sl dual group} of $G$. The dual of $S^1$ turns out to be $\mathbb Z$. Therefore the Fourier transform associates to a function on $S^1$ a function on its dual $\mathbb Z$, the spectrum of $f$.  From now on we shall freely  adopt the convention of referring to concepts related to $S^1$ as {\sl primal} and concepts related to $\mathbb Z$ as {\sl dual} or {\sl spectral}.
\end{remark}

Let us say a few words about about the $L^2$-theory and discuss a little important differences with the $L^1$-theory, within which we will stay in  the rest of this chapter. Define $L^2(S^1,\mathbb C)$ as the set $\{f:S^1\to\mathbb C: \int_0^{2\pi} |f(t)|^2 dt<\infty\}$. $L^2(S^1,\mathbb C)$ is a normed linear space with the norm defined by $\parallel f\parallel_2={1\over 2\pi}( \int_0^{2\pi} |f(t)|^2 dt)^{1\over 2}$. Putting a scalar product on it by mean of $(f,g)= \int_0^{2\pi} f(t)\bar g(t) dt$ makes of $L^2(S^1,\mathbb C)$ a Hilbert space. Functions in $L^2(S^1,\mathbb C)$ will be referred to as {\sl square-summable} (or {\sl square Lebesgue integrable}). Consider  the  dual space  $l^2=\{\alpha:\mathbb Z\to\mathbb C: \sum_{k\in\mathbb Z}|\alpha(k)|^2<\infty\}$. The norm $|\alpha|_2= (\sum_{k\in\mathbb Z}|\alpha(k)|^2)^{1\over 2}$ gives  $l^2$ the structure of a Hilbert space.  Because $S^1$ is of finite measure $L^2(S^1,\mathbb C)\subset L^1(S^1,\mathbb C)$, i.e. square-summable functions on $S^1$ are summable. Therefore the Fourier transform is defined.
The fundamental theorem about the Fourier transform in $L^2$ is the following
\begin{theorem}{\bf (Plancherel's Theorem)} The Fourier transform is an isometric isomorphism between the Hilbert spaces $L^2(S^1,\mathbb C)$ and $l_2$:
\end{theorem}

$$\parallel f\parallel_2=|\hat{f}|_2$$
\begin{flushright}$\square$\end{flushright}

Note that there is no such theorem for $L^1$. In the words of Dym and McKean, quoting word by word from their authoritative treatise (\cite{D-McK}, pp. 42-43):
\bigskip

``The situation for $L^1(S^1)$ is very much more complicated. The information at hand about the class A of Fourier coefficients $\hat f$ of summable functions may be summarized as follows:

\begin{itemize}
\item[ ] (a) \,\,\,A is populated by bounded functions $\hat f$:

$$\parallel\hat f\parallel_\infty\,\,=\,\,\lim_{n\uparrow\infty}(\parallel f^n\parallel_1)^{1/n}\,\,\leq\,\,\parallel f\parallel_1.$$

\item[ ] (b) \,\,\, $\hat f=0$ iff $f=0$, i.e.,\,\,\,$\wedge$\,\, is a 1:1 map.

\item[ ] (c) \,\,\, $\hat f=0$ iff $f=0$\, at $\pm\infty$, i.e., $\lim_{|n|\uparrow\infty}\hat f(n)=0$.
\item[ ] (d )\,\,\, A is an algebra: namely,

$$(\hat f\hat g)(n)=\hat f(n)\hat g(n)= (f\circ g)^{\wedge}(n).$$
\end{itemize}
Unfortunately, (a) and (c) do not suffice to single out precisely which functions $\hat f$ arise as Fourier coefficients of summable functions. The situation is thus entirely different from $L^2(S^1)$ where the condition

$$\parallel\hat f\parallel_2=(\sum|\hat f|^2)^{1/2}<\infty$$

\noindent is decisive. The best information currently available indicates that A does not have {\sl any} neat description.''

\bigskip 

It is therefore quite a surprise that, as we will prove in this chapter, in the particular case of {\sl fuzzy subsets} such a  description is indeed possible. Even if we frankly doubt that  the authors of the book just quoted would consider it {\sl neat}.

\section {Fuzzy subsets of the unit circle}

Let us introduce the main character of this chapter: $\mathcal F=L^1(S^1,\mathcal I)$, where $\mathcal I=[0,1]$. We will refer to elements of $\mathcal F$ as {\sl integrable fuzzy subsets}  of the unit circle. Remark that we are adopting the linguistic convention of considering fuzzy objects a superset of crisp objects. In our case the crisp objects are characteristic functions of subsets of the unit circle. Obviously $\mathcal F$ is no longer a linear space, but it inherits the structure of a metric space with respect to the metric induced by the distance $d(f,g)=\parallel f-g\parallel_1$. That it is actually complete is the content of the following
\begin{theorem} $\mathcal F$ is a complete metric space.
 \end{theorem}
{\it Proof.}  From a Cauchy sequence in $L^1$ we can extract a subsequence which converges almost everywhere to a function $f$, which is obviously still a function with values in $[0,1]$. Because $S^1$ is of finite measure, convergence a.e. implies convergence in measure. Then, being the sequence $\{f_n\}$ trivially dominated by the function identically equal to 1, the theorem of Lebesgue allows  to conclude that $f\in \mathcal F$.   
\begin{flushright}$\square$\end{flushright}

The following theorem gives a necessary condition for a sequence of complex numbers  to be the spectrum of a fuzzy subset.

\begin{theorem} Let $f\in \mathcal F$. Then its Fourier coefficients satisfy the inequalities
\end{theorem}
   $$0\leq c_0\leq 1,\,\,\,|c_k|\leq{\sqrt 2\over\pi},\,\,\,k\in\mathbb Z^*$$ 
{\it Proof}. Being $c_0$ the mean value of $f$, one has trivially that $0\leq c_0\leq 1$.  Let us compute bounds for the Fourier coefficients $a_k,k\in\mathbb N^*$:

 $$a_k(f)={1\over \pi}\int_0^{2\pi} f(t)cosktdt={k\over \pi}\int_0^{2\pi\over k}f(t)cosktdt\leq$$
 $$\leq {k\over \pi}\int_0^{2\pi\over k}f(t)(coskt)^+dt={2k\over \pi}\int_0^{\pi\over 2k}(coskt)^+ dt={2\over\pi}$$
 
 \noindent where we have used the standard notation $g^+$ for the positive part of g. 
 Analogous computations yield the same value $2\over\pi$ as a bound for $b_k$. Therefore  $|c_k|\leq{\sqrt 2\over\pi}$.
 
Let us add the  remark that the bounds are sharp, in the sense that there are functions for which the $k$-th coefficient is equal to the bound: for example, for the signature of the positive part of $coskt$ one has $a_k={2\over\pi}$, and so is the case for the coefficient $b_k$ of the signature of the positive part of $sinkt$.
\begin{flushright}$\square$\end{flushright}

\section{Crisp subsets of order n}
 
Later in the course of the chapter we will find also sufficient conditions for a sequence to be the sequence of the Fourier coefficients of a fuzzy subset. The main ingredient will be a spectral approximation procedure, the idea of which, more than half a century old, goes back to a set of papers by Szasz (\cite{Szasz2}), Friedman (\cite{Fried}) and, especially, Ghizzetti (\cite{Ghizz}). The idea turns around a very special subset of $\mathcal F$, called $\mathcal G$ in honour of Ghizzetti, formed by the characteristic functions of  finite unions of (non degenerated) arcs: a nice class of {\sl crisp} subsets.  To set the notations, an element of $\mathcal G$ is determined by a collection $\Gamma$ of $n$ arcs $\gamma_1,\dots,\gamma_n$ where  $\gamma_r:[0,1]\to S^1,\,\, r=1,\dots,n$. Let us denote the characteristic function  by $\chi_\Gamma$. Of course, it is no restriction to suppose that the left point of $\gamma_1$ is the point $1$ on the unit circle. Using the parametrization of $S^1$ given by $[0,2\pi)$, we can sum up this discussion in a formal 

\begin{definition} 
For every $n\in\mathbb N$, consider the collection of all real numbers such that $0=\xi_1<\eta_1<\xi_2<\eta_2<\dots<\xi_n<\eta_n<2\pi$ and define $[\xi,\eta]=\cup_{r=1,\dots,n}[\xi_r,\eta_r]$. The subset $\mathcal G_n\subset \mathcal F$ is the collection of all the characteristic functions

$$\{\chi_{[\xi,\eta]}:S^1\to\mathcal I\}$$

\noindent Elements of $\mathcal G_n$ are called {\sl crisp subsets of order n}. Elements of $\mathcal G=\cup_{n\in\mathbb N}\,\,\,\mathcal G_n$ are called {\sl crisp subsets of finite order}.
\end{definition}

\noindent The relationship with the description of $\mathcal G$ in terms of arcs is given by $\gamma_r(0)=e^{i\xi_r}$ and $\gamma_r(1)=e^{i\eta_r}$. 
Note that $\mathcal G$ is by no means complete as a metric subspace of $\mathcal F$. It is  very easy  to find  Cauchy sequences which do not  converge to a function of $\mathcal G$:  take a sequence of arcs which collapses to a point; points are {\sl degenerate} arcs, so they do not belong in $\mathcal G$. 

Let us end the section by computing the spectrum $\{c_k\}_{k\in\mathbb Z}$ of a crisp subsets of order $n$, $\chi_{[\xi,\eta]}$:

$$c_0(\chi_{[\xi,\eta]})={1\over{2\pi}}\int_0^{2\pi}\chi_{[\xi,\eta]}(t)dt=\sum_{r=1}^n(\eta_r-\xi_r)$$

$$c_k(\chi_{[\xi,\eta]})={1\over 2\pi}\int_0^{2\pi} \chi_{[\xi,\eta]}(t)e^{ikt}dt={1\over{2\pi ki}}\sum_{r=1}^n(e^{ik\eta_r}-e^{ik\xi_r}),\,\,\, k\in\mathbb Z^*$$

\section{Main theorems, heuristics and structure of the proof}

We are now in a position to state formally one of our two main theorems:

\begin{theorem}
For every $f\in\mathcal F$ there exists a sequence  $\{\chi_n\}_{n\in\mathbb N}, \,\,\,\chi_n\in\mathcal G_n$ such that 

$$c_k(f)=c_k(\chi_n), k=0,\dots,n-1$$
\end{theorem}

\noindent In plain language this amounts to say that to every fuzzy subset of the unit circle we can associate a sequence of crisp subsets of order $n$ having the same first $n$ Fourier coefficients as the original fuzzy subset. The sequence of crisp subsets of order $n$ can be thought of as a progressive reconstruction of the spectrum.

The proof of this theorem will be a byproduct of the proof of our second main theorem, giving a complete characterization of the spectra of fuzzy subsets of the unit circle. Before being able to give a formal statement, we need to define quite a few mathematical concepts, whose introduction we hope to justify, by indulging somehow in some heuristic considerations.

\subsection{Heuristics 1, or how to construct crisp subsets of order $n$ with preassigned  Fourier coefficients}

Given $n$ complex numbers $\alpha_0,\dots,\alpha_{n-1}$, how to find a crisp subset of order $n$, $\chi$ such that $c_k(\chi)=\alpha_k$? First of all, the $\alpha$'s must satisfy the necessary conditions $0\leq\alpha_0\leq 1, |\alpha_k|\leq{\sqrt 2\over\pi}$. It is straightforward to solve the problem for $n=1$: just take the crisp subset of order 1 defined as ${1\over{2\pi}}\chi_{[0,\alpha_0]}$. Remark the simple fact that this is the only solution: characteristic functions of  intervals $[\xi,\xi+\alpha_0]\subset[0,2\pi]$ are not crisp subsets of order 1, according to our definition. Things get immediately much worse when it comes to higher Fourier coefficients. The task of finding a solution $0=\xi_1<\eta_1<\dots\xi_n<\eta_n<2\pi$ to the transcendental equation 

$${1\over{2\pi ki}}\sum_{r=1}^n(e^{ik\eta_r}-e^{ik\xi_r})=\alpha_k$$

\noindent by brute force seems hopeless. One has to find some tricks.

One of the most important ideas, when studying properties of periodic functions, is to extend a functions on $S^1$ to a holomorhic function in the disk $D=\{z\in\mathbb C: |z|<1\}$, hoping that informations on the holomorphic extension can yield informations on the boundary value function, i.e. the original periodic function. This brings us to the concept of Hardy spaces, which we will briefly review in the next section.

\subsection{Hardy spaces}

\begin{definition} The Hardy space $H^1$ is defined as the space of all functions $g$ holomorphic in the unit disc $D=\{z\in\mathbb C: |z|<1\}$ for which the norm

$$\parallel g\parallel_{H^1}=\sup_{r<1}{1\over 2\pi}\int_0^{2\pi}|g(re^{it})|dt$$

\noindent is finite.
\end{definition}

\noindent Note that the radial limit 
$\tilde g(e^{it})=\lim_{r\to 1}g(re^{it})$
 exists almost everywhere, by Fatou's lemma. Moreover $\tilde g\in L^1(S^1,\mathbb C)$ and $\parallel g\parallel_{H^1}=\parallel\tilde g\parallel_1$ (for a proof of these simple facts see, e.g., \cite{Koosis}.) It is also possible to define the space $H^1$ directly as the subspace of the functions in $f\in L^1(S^1,\mathbb C)$ with vanishing negative Fourier coefficients. The function $\tilde g$ with $\tilde g(e^{it})=\sum_{k\in\mathbb N} c_k(g)e^{ikt}$ will be identified with the holomorphic function on $D$ defined by $g(z)=\sum_{k\in\mathbb N} c_k(g)z^k$. Observe  the obvious fact that the correspondence between $L^1(S^1,\mathbb C)$ and $H^1$  defined by
 
 $$\sum_{k\in\mathbb N} c_k(g)e^{ikt}\mapsto\sum_{k\in\mathbb N} c_k(g)z^k$$
 
\noindent is not 1-1: functions having the same non-negative Fourier coefficients, but differing otherwise, have the same image in $H^1$. The situation is different, when one considers $L^1(S^1,\mathbb R)$. In this case the Fourier coefficients verify $c_{-k=}\bar c_k$, so that negative coefficients are determined by positive ones. Therefore the correspondence is 1-1. This allows us to formulate the following

\begin{definition}
The 1-1 map $h:L^1(S^1,\mathbb R)\to H^1$ given by $h(f)(z)=\sum_{k\in\mathbb N} c_k(f)z^k$ is called the {\sl Hardy map}. h(f) is called the {\sl Hardy fuction} associated to $f$.
\end{definition}

The idea then is to consider  Hardy functions associated to  fuzzy subsets and try to see if we can take advantage of the powerful methods of complex analysis to undercover some hidden useful informations. Let us start by computing Hardy functions associated to crisp subsets of order $n$. Let $\chi$ be the characteristic function of the  union of the intervals $[\xi_1,\eta_1],\dots,[\chi_n,\eta_n]$. Remembering that 
$$c_k(\chi)={1\over{2\pi ki}}\sum_{r=1}^n (e^{ik\eta_r}-e^{ik\xi_r})$$

\noindent we have:

$$h(\chi)(z)=\sum_{k=1}^\infty c_k(\chi)z^k=\sum_{k=1}^\infty(-{1\over{2\pi ki}}\sum_{r=1}^n (-e^{ik\eta_r}+e^{ik\xi_r}))z^k=$$

$$=-{1\over{2\pi i}}\sum_{r=1}^n(\sum_{k=1}^\infty (-e^{ik\eta_r}+e^{ik\xi_r})){z^k\over k}=
-\sum_{r=1}^n({1\over{2\pi i}}\sum_{k=1}^\infty \frac{(-ze^{i\eta_r})^k}{k}-\sum_{k=1}^\infty \frac{(-ze^{i\xi_r})^k}{k})=$$

$$=-\sum_{r=1}^n({1\over{2\pi i}}\sum_{k=0}^\infty (-1)^k\frac{(-ze^{i\eta_r})^{k+1}}{k+1}-\sum_{k=0}^\infty 
(-1)^k\frac{(-ze^{i\xi_r})^{k+1}}{k+1})$$

\noindent Recalling that

$$ln(w+1)=\sum_{k=0}^\infty(-1)^k\frac{w^{k+1}}{k+1}$$

\noindent and applying it to, respectively, $w=-ze^{i\eta_r}$ and $w=-ze^{i\xi_r}$, we see that the last expression in our computation is equal to

$$-{1\over{2\pi i}}\sum_{r=1}^n ln(-ze^{i\eta_r}+1)-ln(-ze^{i\xi_r}+1)=
-{1\over{2\pi i}}\sum_{r=1}^n ln\frac{1-ze^{i\eta_r}}{1-ze^{i\xi_r}}.$$

\subsection{Heuristics 2, or what is the price to pay to get rid of logarithms}
Let us rest a moment to comment on what we got. We have an explicit expression relating two different representations of the holomorphic function
$h(\chi)$:

$$\sum_{k=1}^\infty c_k(\chi)z^k= -{1\over{2\pi i}}\sum_{r=1}^n ln\frac{1-ze^{i\eta_r}}{1-ze^{i\xi_r}}$$

\noindent In the left-hand  we have the Fourier coefficients of $\chi$, in the right-hand the end-points of the arcs defining $\chi$. The problem is the logarithm, which renders the relation transcendental, whatever that means. Surely not friendly to the perspective of actually explicitly computing  the arcs as functions of the Fourier coefficients. So let us get rid of it. As it is known to one and all, getting rid of logarithms is a children's game: it is enough to exponentiate. Let us do that thoughtlessly. We will figure out later how to pay the bill. While we are at that, let us also get rid of the factor ${1\over{2\pi i}}$:

$$\mbox{exp}(-2\pi i(\chi)(z))=\mbox{exp}(\sum_{r=1}^n ln\frac{1-ze^{i\eta_r}}{1-ze^{i\xi_r}})=\prod_{r=1}^n\mbox{exp}( ln\frac{1-ze^{i\eta_r}}{1-ze^{i\xi_r}})=\prod_{r=1}^n\frac{1-ze^{i\eta_r}}{1-ze^{i\xi_r}}$$

\noindent Let us rest again and comment: now the right-hand side is really beautiful. It is a rational function, the quotient of two polynomials of degree $n$, holomorphic on the disk, the poles being on the conjugates of the end points of the arcs defining the crisp set $\chi$. A little problem on the left-hand: no longer the Hardy function, but still a holomorphic function on the disk. Its Taylor coefficients are no longer the Fourier coefficients of $\chi$. In exchange their relationship with the geometric data of the arcs defining $\chi$ seems very promising. If we call $s_k$ the Taylor coefficients of the function $e^{-2\pi i}\chi$, we have that

$$\sum_{k\in\mathbb N}s_k z^k= \prod_{r=1}^n\frac{1-ze^{i\eta_r}}{1-ze^{i\xi_r}}$$

\noindent Having found such a nice formula, the problem is shifted to the study of the relationship between  the Fourier coefficients  of a function $f\in L^1(S^1,\mathbb R)$ and  the Taylor coefficients of the holomorphic function $e^{-2\pi h(f)}$. In the next section we will pay the bill.

\subsection{Nonlinear Fourier transforms and nonlinear coefficients}
Take a function $f\in L^1(S^1,\mathbb R)$ and associate to it the holomorphic function on the disk $e^{-2\pi h(f)}$. Let the development of $e^{-2\pi h(f)}$ in power series be

$$e^{-2\pi h(f)}= \sum_{k=0}^{\infty}a_kz^k\]

\noindent Remark that the $0-th$ term of the series is $1$. Define $s_0=2\sin\pi c_0$ and rescale the coefficients according to  $a_k=\frac{-s_k}{ie^{\pi ic_0}}$. One has that

$$e^{-2\pi h(f)}=1-i e^{\pi ic_0}\sum_{k=1}^{\infty}s_kz^k$$

\noindent We are ready to state the following
\begin{theorem} The transformation rules between the coefficients $s_k$ and $c_k$ are algebraic. Precisely $s_k=P(c_1,\dots,c_{k-1})$ and $c_k=Q(s_1,\dots,s_{k-1})$, with $P$  and $Q$ polynomials of degree $k$.
\end{theorem}

{\sl Proof.} We will actually prove that 

$$P(c_1,\dots,c_{k-1})=2\pi[e^{-\pi ic_0}c_k-i\sum_{r=1}^{k-1}(1-\frac{r}{k})c_{k-r}s_r]$$ 

\noindent The equation involving $Q$ is analogous. Taking the logarithm of  both sides of the identity

\[e^{-2\pi ih_f(z)}=1-i e^{\pi ic_0}\sum_{k=1}^{\infty}s_kz^k\]
 
\noindent we obtain
\[-2\pi i\sum_{k=1}^{\infty}c_kz^k=\ln[1-i e^{\pi ic_0}\sum_{k=1}^{\infty}s_kz^k]\]
and, by  differentiating on both sides
\[-2\pi i(c_1+\sum_{k=2}^{\infty}c_kkz^{k-1})=\frac{-i e^{\pi ic_0}(s_1+\sum_{k=2}^{\infty}s_kkz^{k-1})} {1-i e^{\pi ic_0}\sum_{k=1}^{\infty}s_kz^k}\]
Let us first factor out $z^k$ for all $k\geq 0$:
\[\sum_{k=0}^{\infty}\Omega_kz^k=0\]
We obtain, at order $0$ 
\[s_1 e^{-\pi ic_0}=2\pi c_1\]
and at order $k-1,\ k>1$ 
\[s_k=2\pi[e^{-\pi ic_0}c_k-i\sum_{r-1}^{k-1}(1-\frac{r}{k})c_{k-r}s_r]\]
\begin{flushright}$\square$\end{flushright}

\begin{definition} Given $f\in L^1(S^1,\mathbb R)$ we  call the Taylor coefficients of $e^{-2\pi h(f)}$ the {\sl nonlinear coefficients} of $f$.
\end{definition}

Let us call $\mathcal N$ the map which associates to a function $f$ its nonlinear coefficients. One could think of this map as a {\sl nonlinear Fourier transform} and of the set of nonlinear coefficients as a {\sl nonlinear spectrum}. The effect of this nonlinearization on the original problem is to simplify the relationship between coefficients and geometric data of crisp subsets of order $n$. Therefore the strategy will be to obtain a characterization of the nonlinear coefficients of fuzzy subsets and then transfer the conclusions to the Fourier coefficients, by mean of the algebraic maps which relate the ones to the others. The abstract scheme is illustrated by the following diagram

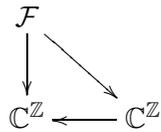
\begin{displaymath}
\xymatrix{
\mathcal F \ar[d] \ar[dr]  &  \\
\mathbb C^\mathbb Z & \mathbb C^\mathbb Z\ar[l] }
\end{displaymath}

\noindent where the vertical arrow is the Fourier transform, the diagonal arrow the nonlinear Fourier transform and the horizontal arrow the algebraic map between  nonlinear and Fourier coefficients.

We are now ready to introduce the last main ingredient of our proof and finally give a statement of our main theorem.

\subsection{Toeplitz matrices}

\begin{definition} A square matrix  $A=(a_{ij})$ is called a Toeplitz matrix if the entries $a_{ij}$ depend only on $i-j$.
\end{definition}

\noindent Note that this is equivalent to saying that the entries are constant on the NW to SE diagonals.

There are (infinitely) many ways to associate to a sequence of complex numbers $s:\mathbb Z\to\mathbb C$ a Toeplitz matrix. We will restrict ourselves to the hermitian case ($s_{-k}=\bar s_k$) and define two sequences $T_k=T_k(s)$ and $W_k=W_k(s),\,\,k\in\mathbb N$ of, respectively $(k+1)\times(k+1)$ and $k\times k$  Toeplitz matrices  by

 \[T_k(s)=\left(\begin{array}{cccc}s_0 & s_1 & s_2 & \dots s_k \\s_{-1} & s_0 & s_1 & \dots s_{k-1}  \\\dots & \dots & \dots & \dots \\s_{-k} & s_{-k+1} & s_{-k+2} & s_0\end{array}\right)\]
 
 \noindent
 
 \[W_k(s)=\left(\begin{array}{cccc}s_1 & s_2 & s_3 & \dots s_k \\s_0 & s_1 & s_2 & \dots s_{k-1}  \\\dots & \dots & \dots & \dots \\s_{-k+2} & s_{-k+3} & s_{-k+4} & s_1\end{array}\right)\]

\noindent Remark that $W_k$ is obtained from $T_k$ by suppressing the first column and the last line and that $T_k$ is hermitian, while $W_k$ is not. Call $D_k$ the determinant of $T_k$ (which is real, being $T_k$ hermitian) and $F_k$ the determinant of $W_k$. We state three important results about the Toeplitz sequences. For the missing proofs, we refer to the historical papers of Szeg\H o (\cite{Szeg1} and \cite{Szeg2}) or to the more recent book \cite{Heinig-Rost}.

\begin{lemma}
$D_{k-1}^2-D_{k-2}D_k=|F_k|^2$
\end{lemma}
\begin{flushright}$\square$\end{flushright}

\begin{lemma}
Suppose that $D_0>0,\,\,D_1>0\,\,,\dots,D_{n-1}>0$ and that $D_k=0\,\,\forall k\geq n$,
then the zeroes $\alpha_1,\dots,\alpha_n$ of the following polynomial $P(z)$ of degree $n$

\[P(z)=\left|\begin{array}{cccccc}1 & z & z^2 & \dots & z^{n-1} & z^n \\s_0 & s_1 & s_2 & \dots & s_{n-1} & s_n \\s_{-1} & s_0 & s_1 & \dots & s_{n-2} & s_{n-1} \\\dots & \dots & \dots & \dots & \dots & \dots \\s_{-n+1} & s_{-n+2} & s_{-n+3} & \dots & s_0 & s_1\end{array}\right|\]
are all simple and of modulus 1. Moreover $s_k$ can be written in a unique way as 

\[s_k=\sum_{r=1}^n\mu_r\alpha_r^k,\,\,\,k\in\{-n+1,-n+2,\dots, n-2,n-1\}\]

\noindent where the numbers $\mu_1,\dots,\mu_n$ are real and positive.
\end{lemma}
\begin{flushright}$\square$\end{flushright}

\subsection{Heuristics 3, or what the Toeplitz matrices and determinants suggest for the solution of our problem}
Imagine that the sequence $s_k$ of the nonlinear coefficients of a crisp set of order $n$ satisfy the condition of Lemma 2, i.e. all the Toeplitz determinants up to order $n-1$ are (strictly) positive and definitively zero from order $n$, then the roots of the algebraic equation $D_n=0$ are $n$ distinct points on the unit circle, call them $e^{it_1},\dots,e^{it_n}$ and you are able to write the first $n$ nonlinear coefficients as linear combinations (with positive coefficients) of rotations $e^{ikt_r}$ of the  $n$ roots. Very nice indeed! That this is actually the case will be the main ingredient of the proof of our main theorem, which we are, at last, ready to state, after one more definition.

\begin{definition} Let $\mathbb C^\mathbb Z_H$ be the space of all hermitian sequences of complex numbers. We say that a sequence $\beta\in\mathbb C^\mathbb Z_H$ is of {\sl order} $n$, if the sequence of the Toeplitz determinants $D_k(\beta)$ verifies

$$D_k(\beta)>0,\,\,\,k=0,\dots,n-1,\,\,\,D_k(\beta)=0,\,\,\,\forall k\geq n$$

\noindent We denote by $V_n$ the subset of all  sequences of order $n$. Define $V=\cup_{n=1}^\infty\,\,V_n$. Elements of $V$ are called the sequences of {\sl finite order}. 
Define now the set $V_\infty$ of all sequences $\beta\in\mathbb C^\mathbb Z_H$ which verify 

$$D_k(\beta)>0\,\,\,\forall k\in\mathbb N$$

\noindent We call  elements of $V_\infty$  sequences of {\sl infinite order}. Finally we denote by $V$ the union of the sets of sequences of finite order and of infinite order:

$$V=(\cup_{n=1}^\infty\,\,V_n)\,\cup\,V_\infty$$

\end{definition}

We are now in a position to state our main theorem:

\begin{theorem} 
$\mathcal N(\mathcal F)=V$ and $\,\,\mathcal N(\mathcal G_n)=V_n$
\end{theorem}

Translated in plain language, the theorem says that the sequence of nonlinear coefficients of a fuzzy set either is of finite order, and in this case the fuzzy set is actually a crisp set of order $n$, or it is of infinite order. Stated otherwise, the nonlinear Fourier transform maps $\mathcal F$ into $V$. Viceversa, every sequence in $V$ arises as the sequence of the nonlinear coefficients of a fuzzy set $f\in\mathcal F$, so that the map $\mathcal N$ is actually a bijection. Remark that the proof that we will give of the surjectivity is constructive, in the sense that we will be able to build a sequence of crisp sets having the nonlinear coefficients equal up to order $n$ to the nonlinear coefficients of a preassigned fuzzy set. Therefore we will get at the same time a proof of the spectral approximation theorem stated at the beginning of Section 4.

\section {Proof of the main theorems}

For the readers convenience we start by highlighting two main facts, somehow hidden in the flow of computations, before starting the actual proof.

\begin{itemize}
\item $\mathcal N(\mathcal G_n)=V_n$. Moreover, if $f\in\mathcal G_n$ is the characteristic function of the segments $[\xi_r,\eta_r],\,\,r=1,\dots,n$ and $s_k$ its nonlinear coefficients, the roots of the Toeplitz polynomial $D_n(s)$ are $e^{i\xi_1},\dots,e^{i\xi_n}$ and $W_n(s)=D_{n-1}(s) e^{i\sum_{r=1}^n\xi_r}$

\item $\mathcal N(\mathcal F\backslash \mathcal G)=V_\infty$
\end{itemize}

\noindent{\sl Proof of the direct part}: $\mathcal{N}(\mathcal{F})\subset V$

\bigskip

Let $f\in\mathcal{F}$ and, as usual, let  $c_k$ be the sequence of its Fourier coefficients. Let us first notice that
\[c_0=0\ (or\ c_0=1) \Rightarrow  f\equiv 0 \ a.e.\Rightarrow c_1=c_2=\dots=0\]

\noindent This implies that 
\[s_0=s_1=\dots=0\Rightarrow D_0=D_1=\dots=0\]
hence $\mathcal{N}(f)\in V_1$. Let us suppose that 

\[0<c_0<1\]
and 
\[D_0>0,D_1>0,\dots,D_{n-1}>0.\]

Set $\mathcal{G}_n(\lambda)=\{\chi(\xi_1,\xi_2,\dots,\xi_n)\in\mathcal{G}_n:\sum_{i=1}^n\xi_i=\lambda\}$. We are going to show that $D_n\geq 0$,  the equality holding if and only if $f\in\mathcal{G}_n(arg(W_n))$ and its first $n$ Fourier coefficients are equal to $c_0,c_1,\dots,c_{n-1}$. This is equivalent to say that  $\mathcal{N}(f)\in V$, and that $\mathcal{N}(f)\in V_n\,$ if and only if $D_n=0$.

 \begin{remark}
In plain words this means that if $D_n$ vanishes (and consequently, by the previous observation,  $D_k=0,\ k>n$ ) then $f$ must be a crisp set of order $n$, while, if $D_n\neq 0$, then $f$ does not coincide with any crisp set up to order $n$, so we can look at the next $Toeplitz$ determinant $D_{n+1}$, and so on. In our notation, this is equivalent to say that $f\in\mathcal{N}(\mathcal{F}\backslash\mathcal{G})$.

\end{remark}

Consider an arbitrary $\lambda\in\mathbb{R}$. We can claim that there exists $\chi(\xi_1,\xi_2,\dots,\xi_n)\in\mathcal{G}_n(\lambda)$ such that 
its first $n$ Fourier coefficients are equal to $c_0,c_1,\dots,c_{n-1}$ and 
\[\sum_{i=1}^n\xi_i\equiv\lambda(mod\ 2\pi)\]

\noindent From this we obtain that 
\[\sum_{r=1}^n(\eta_r+\xi_r)=2\pi c_0+2\lambda\]
Define 
\[P(x)=2^{2n}\prod_{r=1}^nsen\frac{x-\xi_r}{2}=\sum_{k=-n}^{n}P_ke^{ikx},\ (P_{-k}=\overline{P}_k)\]
Since $P_n=(-1)^ne^{-\frac{1}{2}i\sum_{i=1}^n(\xi_r+\eta_r)}=(-1)^ne^{-i(\pi c_0+\lambda)}$, then 
\[P(x)>0\mbox{ on } (\xi_i,\eta_i)\]
\[P(x)<0\mbox{ on } (\eta_i,\xi_{i+1})\]
We notice that the fact that $\lambda$ is defined $mod\ 2\pi$ does not affect the above identities since $\lambda$ appears as an exponent.

\noindent The construction of the trigonometric polynomial $P$ with the above property is crucial to our analysis: we notice that the function $[f-\chi_n(\lambda)]P(x)$ coincides with $sign[P]$ over the period $(\xi_1,\xi_1+2\pi)$.\\
Therefore we have 

\[\frac{1}{2\pi}\int_{\xi_1}^{\xi_1+2\pi}[f(x)-\chi_n(\lambda)]P(x)dx\geq0\]
where the equality holds if and only if $f= \chi_n(\lambda)$.

The {\sl magic} produced by $P$ consists in the fact that, when multiplied by $f$, we get an expression where the Fourier coefficients of $f$ appear, and, thanks to the boundedness of $f$ we are able to obtain inequalities. More precisely, by expanding $P$, the above inequality, after some involved, even though basic, algebraic manipulations, reduces to

\[D_{n-1}\geq Re(e^{-i\lambda}F_n)\]
where the equality holds if and only if $f= \chi(\lambda)$.\\
Since neither $D_{n-1}$ nor $F_n$ depend on $\lambda$, the above equation can be written as
\[{n-1}\geq|F_n| cos(\omega-\lambda)\]

\noindent where $arg(F_n)=\omega$. From which we deduce  ($\lambda$ is arbitrary) that

\[D_{n-1}\geq F_n\]
If $D_{n-1}>|F_n|\ \forall \lambda$, then $f$ does not belong to any $\mathcal{G}_n(\lambda)$.\\
On the other hand, if $D_{n-1}=F_n$, then (4) gives an equality only when $\lambda=\omega$, in which case $f$ coincides with $\chi(\lambda)$.\\
Remembering that
\[D_{n-1}^2-D_{n-1}D_n =|F_n|^2\]
and that
\[D_{n-2}>0\]
we can conclude that 
\[D_n\geq 0\]
where the equality holds if and only if $f= \chi(\lambda)$.\\
\bigskip

\noindent{\sl Proof of the inverse part}: $V\subset\mathcal{N}(\mathcal{F})$
\bigskip

We already know, from the direct part of the proof, that $\mathcal N(\mathcal G)=\cup_{n=1}^\infty V_n$. Therefore let us pick a sequence $\beta\in V_\infty$. Our goal is to construct $f\in\mathcal{F}$ such that $\mathcal{N}(f)=\beta$.
We know that, for any $\lambda$, there exists $\chi_n(\lambda)\in\mathcal{G}_n$ such that 
\[\frac{1}{2\pi }\int_0^{2\pi} \chi_n(x,\lambda)e^{ikx} dx=c_k,\ \ (k=0,1,\dots,n-1)\]
Make the choice  $\lambda=0$ once for all and allow $\lambda$ to disappear from the notation.
Consider the sequence $\chi_k$ and define, on $[0,2\pi]$, the function 
\[\Theta_n(x)=\int_0^x\chi_n(t)dt\]
We observe that the $\Theta_n$ are nonnegative and nondecreasing $\forall n$ and that $\Theta_n$ is uniformly bounded, since
\[0\leq\Theta_n(x)\leq\Theta_n(2\pi)=\int_0^{2\pi}\chi_n(t)dt
=2\pi c_0\]
and  uniformly Lipschitz, since
\[|\Theta_n(x'')-\Theta_n(x')|=|\int_{x'}^{x''}\chi_n(t)dt|\leq|x''-x'|\]
By Ascoli-Arzel\`a, we can extract from  $\Theta_n$ a subsequence $\Theta_{n_h}$
converging to some function $F$, nonnegative, nondecreasing and Lipschitz, since 
\[|F(x'')-F(x')|\leq|x''-x'|\]
Therefore $F$ is differentiable almost everywhere. Define
\[f=F'\]
Note that  $0\leq f\leq 1$
Having constructed a fuzzy set $f\in\mathcal F$, it remains to show that $\mathcal{N}(f)\in V$. Since, by construction,  $\mathcal{N}[\chi(x,\lambda)]\in V_n$, we only need to show  that  the $k$-th nonlinear coefficient of $f$ coincides with the $k$-th nonlinear coefficient of $\chi(x,k)$. In other words, what is left to be proved is that

\[\frac{1}{2\pi}\int_0^{2\pi}e^{ikx}dF(x)=c_k,\]
Let us  notice that, from the first part of the theorem, we have  that definitively
\[\frac{1}{2\pi}\int_0^{2\pi}e^{ikx}d\Theta_n(x)=c_k\]
This implies that
\[\lim_{n\to\infty}\frac{1}{2\pi}\int_0^{2\pi}e^{ikx}d\Theta_n(x)=c_k\]
and, passing to the subsequence: 
\[\lim_{h\to\infty}\frac{1}{2\pi}\int_0^{2\pi}e^{ikx}d\Theta_{n_h}(x)=c_k\]

\noindent Integrating by parts 
\[\int_0^{2\pi}e^{ikx}d[\Theta_{n_h}(x)-F(x)]=-ik \int_0^{2\pi}e^{ikx}[\Theta_{n_h}-F(x)]dx \]
because $\forall h\ \  \Theta_{n_h}(0)=0,\ \ \Theta_{n_h}(2\pi)=2\pi c_0$, and then
$F(0)=0$ and $F(2\pi)=2\pi c_0$.\\
We can pass to the limit $h\rightarrow\infty$ under the integral sign, because 
\[|e^{ikx}[\Theta_{n_h}(x)-F(x)]|=|\int_0^x\chi_{n_h}(x)-f(x)dx|\leq\int_0^x dx\leq 2\pi\]
Since, for $h\rightarrow\infty$, $\Theta_{n_h}(x)\rightarrow F(x)$, 
we conclude that 
\[\lim_{h\rightarrow 0}\int_0^{2\pi}[e^{ikx}(\Theta_{n_h}(x)-F(x)]dx=0\]

\begin{flushright}
$\square$
\end{flushright}

\section{On the range of validity of the theorem}

The theorem just proved is only valid {\sl verbatim} for summable fuzzy subsets of the unit circle. Nonetheless, let us remark that every summable fuzzy subset of the real line with compact support, can be extended by periodicity and therefore be interpreted as a fuzzy subset of the unit circle. But something can be said also for arbitrary fuzzy subsets of $\mathbb R$ in  the Schwartz class of smooth functions rapidly decreasing at infinity, having recourse to one of the most celebrated {\sl tricks} of harmonic analysis: Poisson periodization. Let us recall the fundamental facts about it (for the framework and the missing proofs we refer, as usual, to \cite{D-McK}.)

\begin{definition} The {\sl Schwartz class} $\mathcal S$ is the class of all functions $f:\mathbb R\to\mathbb R$ infinitely differentiable and such that

$$\lim_{|x|\to\infty}x^pf^{(q)}(x)=0,\,\,\,\forall p\in\mathbb N, \forall q\in\mathbb N$$

\end{definition}

\begin{definition} Let $f$ be a function in the Schwartz class. Then the series $$\sum_{k\in\mathbb Z} f(x+2\pi k)$$ 

\noindent converges uniformly for $0\leq x<2\pi$ to a smooth $2\pi$-periodic function $f^0$, called the {\sl Poisson periodization} of $f$.

\end{definition}

The fundamental theorem about the Poisson periodization is the following

\begin{theorem} {\bf (Poisson's summation formula)} 

$$f^0(x)={1\over2\pi}\sum_{k\in\mathbb Z}\hat f(n)e^{ikx}$$

\end{theorem}

\noindent{\sl Proof}. Just compute the Fourier coefficients of $f^0$:

$$c_k(f^0)={1\over2\pi}\int_0^{2\pi}f^0(x)e^{ikx}dx={1\over2\pi}\sum_{k\in\mathbb Z}\int_0^{2\pi}f(x+2\pi k)e^{ikx}dx=$$

$$={1\over2\pi}\sum_{k\in\mathbb Z}\int_{2\pi k}^{2\pi(k-1)}f(x)e^{ikx}dx={1\over2\pi}\int_{-\infty}^{+\infty}f(x)e^{ikx}dx={1\over2\pi}\hat f(k)$$

\begin{flushright}
$\square$
\end{flushright}

\noindent The global functioning of the whole {\sl Poisson machinery} can be expressively 
represented by the following commutative diagram:

\begin{displaymath}
\xymatrix{ f \ar[r]\ar[d] & f^0 \ar[d] \\
\hat f \ar[r] & \hat f^0 }
\end{displaymath}

\noindent where the vertical arrows represent   the Fourier transforms of, respectively,  $f$ (a function in the Schwartz class) and  $f^0$ (a converging sequence in $l^1$), while the horizontal arrows represent, respectively, the Poisson periodization and the {\sl sampling operator} on $\mathbb Z$. Therefore, an application of our theorem to the present situation gives the following

\begin{theorem} If $f$ is a fuzzy subset of $\mathbb R$ and $f\in\mathcal S$ then the sampling over $\mathbb Z$ of its spectrum is completely determined by a sequence of spectra of crisp subsets of order $n$ of the unit circle.
\end{theorem}
\begin{flushright}
$\square$
\end{flushright}

We are confident that this theorem could play a role both in the theoretical aspects and in the applications to engineering of defuzzification procedures.

There are other possibilities of extending our main  theorem beyond the original framework of periodic functions. As we mentioned before, Fourier analysis works very well on locally compact abelian groups. Significant parts of the theoretical framework that we used in our proof can be exported to some locally compact abelian groups as well. For example, Hardy spaces on toruses have been recently the object of increasing  investigations, not to mention the great amount of researches on Toeplitz operators arising in different geometries. For some first results in this direction, we refer to the forthcoming paper \cite{cafcat}.

\end{document}